\documentclass[12pt]{article}
\usepackage{amsmath,amssymb,amsthm}
\usepackage{mathtools}
\usepackage[margin=1in]{geometry}
\usepackage{hyperref}
\usepackage{tikz}
\usepackage{graphicx}
\usepackage{url}
\usepackage[mathlines]{lineno}
\usepackage{multirow}
\usepackage{dsfont} 
\usepackage{color}
\usepackage{multicol}
\definecolor{red}{rgb}{1,0,0}

\newtheorem{thm}{Theorem}[section]
\newtheorem{cor}[thm]{Corollary}
\newtheorem{lem}[thm]{Lemma}
\newtheorem{prop}[thm]{Proposition}
\newtheorem{conj}[thm]{Conjecture}
\newtheorem{obs}[thm]{Observation}

\newtheorem{rem}[thm]{Remark}

\theoremstyle{definition}

\theoremstyle{example}
\newtheorem{ex}[thm]{Example}

\def\mtx#1{\begin{bmatrix} #1 \end{bmatrix}}

\DeclareMathOperator{\epr}{epr}

\DeclareMathOperator{\rank}{rank}

\newcommand{\R}{\mathbb{R}}
\newcommand{\F}{\mathbb{F}}

\newcommand{\Rnn}{\R^{n\times n}}

\newcommand{\ZZ}{\mathbb{Z}_2}

\newcommand{\bit}{\begin{itemize}}
\newcommand{\eit}{\end{itemize}}
\newcommand{\ben}{\begin{enumerate}}
\newcommand{\een}{\end{enumerate}}
\newcommand{\beq}{\begin{equation}}
\newcommand{\eeq}{\end{equation}}
\newcommand{\bea}{\begin{eqnarray*}}
\newcommand{\eea}{\end{eqnarray*}}
\newcommand{\bpf}{\begin{proof}}
\newcommand{\epf}{\end{proof}\ms}
\newcommand{\bmt}{\begin{bmatrix}}
\newcommand{\emt}{\end{bmatrix}}
\newcommand{\ms}{\medskip}

\newcommand{\ba}{\begin{array}}
\newcommand{\ea}{\end{array}}
\newcommand{\OL}{\overline}

\usepackage{enumitem}

\newcommand{\zonn}{\{0,1\}^{n \times n}}
\newcommand{\alphann}{\{-\alpha, \alpha\}^{n \times n}}

\begin{document}
\title{On 0--1 matrices whose inverses have entries of the same modulus}
\author{Xavier Mart\'inez-Rivera\thanks{Department of Mathematics and Statistics, University of Victoria,  Victoria, BC, V8W 2Y2, Canada  (martinez.rivera.xavier@gmail.com).
}
}


\maketitle

\begin{abstract}
A conjecture of Barrett, Butler and Hall may be stated as follows:
If $n \geq 3$ and $A \in \{0,1\}^{n \times n}$
(the family of $n \times n$ \ 0--1 matrices)
is a nonsingular symmetric matrix,
then the following two statements are equivalent:
(a) All of the principal minors of $A$ of order $n-2$ are zero;
and 
(b) $A^{-1}$ is a matrix
all of whose entries have the same modulus and
all of whose diagonal entries are equal.
We show that this conjecture holds if $A$ does not have
both a zero and a nonzero principal minor of order $n-4$
(if $n \geq 5$).
The parity of the principal minors of
nonsingular symmetric matrices
$A \in \{0,1\}^{n \times n}$ whose
principal minors of order $n-2$ are all zero is explored,
establishing, in particular,
that the determinants of such matrices are all even.
The aforementioned conjecture was stated in terms of
the enhanced principal rank characteristic sequence
(epr-sequence) of the matrix $A$, a sequence defined as follows:
$\ell_1 \ell_2 \cdots \ell_n$, where $\ell_j$ is either
{\tt A}, {\tt S}, or {\tt N}, based on whether
all, some but not all,  or none of its
principal minors of order $j$ are nonzero.
The epr-sequences that
start with {\tt AN} but whose fourth letter is not $\tt S$
(which are of particular relevance for us)
are completely characterized;
moreover, it is shown that,
if $B$ is an arbitrary symmetric matrix with zero diagonal,
then the odd-girth of its graph,
as well as whether or not its graph is bipartite or an odd cycle,
is deducible from its epr-sequence.
For an arbitrary (not necessarily symmetric)
nonsingular matrix $A \in \zonn$ with $n\geq 3$,
we establish necessary conditions for $A^{-1}$ to be a
matrix all of whose entries have the same modulus;
examples of such conditions are the following:
each row and column of $A$ has an 
even number of nonzero entries;
each entry of $A^{-1}$ is the reciprocal of an even integer;
$\det(A)$ is even;
the difference between any two rows of $A$, 
as well as the difference between any two columns of $A$, 
has an even number of nonzero entries;
if $A$ is symmetric, then $A$ has an 
even number of nonzero diagonal entries;
if $A$ is symmetric and 
$\vec{a}_k$ is the $k$th column of $A$, 
then $A-\vec{a}_k\vec{a}_k^T$
has an even number of nonzero diagonal entries. 
\end{abstract}


\medskip


\section{Introduction}\label{s:intro}
$\null$
\indent
We are concerned with  
nonsingular 0--1 matrices whose inverses are 
matrices all of whose entries have the same modulus
(absolute value), with our chief concern being a 
conjecture of Barrett, Butler and Hall \cite{BBH}. 
The conjecture states that, 
if $n \geq 3$ and 
$A \in \Rnn$ is a nonsingular symmetric 0--1 matrix, 
then the following two statements are equivalent:
(i) All of the principal minors of $A$ of order $n-2$ are zero;
and
(ii) $A^{-1}$ is a matrix
all of whose entries have the same modulus and
all of whose diagonal entries are equal.
Although our primary motivation is 
the aforementioned conjecture, 
which is concerned with symmetric matrices, 
some of our results apply to non-symmetric matrices 
(see Section \ref{s:equim}).

In the interest of keeping our presentation consistent with 
that of \cite{BBH}, and because 
the above formulation of the conjecture lacks simplicity, 
we shall restate the aforementioned conjecture using the terminology used in \cite{BBH},
after introducing the necessary terminology.


For a given positive integer $n$, $[n]:=\{1,2, \dots,n\}$.
For  $B \in \Rnn$ and $\mu, \omega \subseteq [n]$,
$B[\mu, \omega]$ denotes the submatrix of $B$ lying in
rows indexed by $\mu$ and columns indexed by $\omega$,
and $B(\mu, \omega)$ denotes the submatrix obtained by
deleting the rows indexed by $\mu$ and columns indexed by
$\omega$;
$B[\mu, \mu]:= B[\mu]$ and is called a
\textit{principal} submatrix;
and $B(\mu, \mu):= B(\mu)$.
We say that an $n \times n$ matrix has {\em order} $n$.
The determinant of a $k \times k$ principal submatrix is a
\textit{principal} minor, and such a minor has \textit{order} $k$.
The rank of a symmetric matrix is called \textit{principal}
(because of a well-known fact, which is encapsulated in
Theorem \ref{thm: rank of a symm mtx}).

For a given symmetric matrix $B \in \Rnn$,
the \textit{enhanced principal rank characteristic sequence}
(\textit{epr-sequence}) of $B$
(which was introduced in \cite{EPR} and motivated by \cite{P}) is $\epr(B) = \ell_1\ell_2 \cdots \ell_n$, where
\begin{equation*}
   \ell_j =
    \begin{cases}
             \tt{A} &\text{if all of the principal minors of order $j$ are nonzero;}\\
             \tt{S} &\text{if some but not all of the principal minors of order $j$ are nonzero;}\\
             \tt{N} &\text{if none of the principal minors of order $j$ are nonzero (i.e., all are zero);}
         \end{cases}
\end{equation*}
for convenience, $[\epr(B)]_j := \ell_j$, for all $j \in [n]$.
Following the convention adopted in \cite{BBH},
for all $k \in [n]$, we say that
``$B$ starts with $\ell_1\ell_2 \cdots \ell_k$,'' and that
``$B$ ends with $\ell_k\ell_{k+1} \cdots \ell_n$.''
For a given sequence
$t_{i_1}t_{i_2} \cdots t_{i_k}$ from $\{\tt A,N,S\}$,
$\overline{t_{i_1}t_{i_2} \cdots t_{i_k}}$ indicates that
the sequence may be repeated as many times as desired
(or it may be omitted entirely).

The set containing all of the matrices in $\R^{m \times n}$
all of whose entries are either
$\alpha_1, \alpha_2, \dots, \alpha_{k-1}$ or $\alpha_k$,
for some $\alpha_1, \alpha_2, \dots, \alpha_{k} \in \R$,
is denoted by
$\{\alpha_1, \alpha_2 \dots, \alpha_k\}^{m \times n}$.
In particular, matrices in $\{0,1\}^{m \times n}$ are
referred to as 0--1 matrices.
Following the convention adopted in \cite{BBH},
we call a matrix $B \in \Rnn$  \textit{equimodular} if
all of its entries are of the same modulus;
that is, $B$ is equimodular if, for some $\alpha \in \R$,
$B \in \alphann$.
If the diagonal entries of a given square matrix are all equal,
then we say that the matrix has ``constant diagonal.''

As no real symmetric matrix has 
$\tt NNA$ or $\tt NSA$ as a subsequence 
(see \cite[Theorem 2.3]{EPR} and \cite[Corollary 2.7]{EPR}),
the following observation was made in \cite{BBH}:

\begin{obs}\label{n-2 and NAA}
{\rm \cite[p. 42]{BBH}}
Let $n \geq 3$ be an integer and
$A \in \zonn$ be symmetric and nonsingular.
Then $A$ ends with $\tt NAA$
if and only if
all of the principal minors of $A$ of order $n-2$ are zero. 
\end{obs}

Observation \ref{n-2 and NAA} permits an 
alternative formulation of the aforementioned conjecture:

\begin{conj}\label{NAA Conjecture}
{\rm \cite[p. 43]{BBH}}
Let $n \geq 3$ be an integer and
$A \in \zonn$ be symmetric and nonsingular.
Then the following statements are equivalent:
\ben
\item[$(a)$] $A$ ends with $\tt NAA$.
\vspace{-3mm}
\item[$(b)$] $A^{-1}$ is equimodular with constant diagonal.
\een
\end{conj}

For the nonsingular symmetric 0--1 matrix $A$ below,
which is juxtaposed with its inverse, $A^{-1}$,
the following is true:
$A$ ends with $\tt NAA$ ($\epr(A) = \tt NSNAA$) and
$A^{-1}$ is equimodular with constant diagonal
($A^{-1} \in \big\{-\frac{1}{2},\frac{1}{2}\big\}^{5 \times 5}$  and each diagonal entry of $A^{-1}$ is equal to $\frac{1}{2}$):
\begin{equation}
A=
\begin{pmatrix*}[r]
 0 & 1 & 0 & 0 & 1 \\
 1 & 0 & 1 & 0 & 0 \\
 0 & 1 & 0 & 1 & 0 \\
 0 & 0 & 1 & 0 & 1 \\
 1 & 0 & 0 & 1 & 0
\end{pmatrix*};\label{5x5-matrix-example}
\quad
A^{-1} =
\frac{1}{2}
\begin{pmatrix*}[r]
 1 & 1 & -1 & -1 & 1 \\
 1 & 1 & 1 & -1 & -1 \\
 -1 & 1 & 1 & 1 & -1 \\
 -1 & -1 & 1 & 1 & 1 \\
 1 & -1 & -1 & 1 & 1
\end{pmatrix*}.
\end{equation}

By means of an exhaustive computer search,
it was determined by the authors of \cite{BBH},
and announced in \cite[p. 53--54]{BBH}, that
Conjecture \ref{NAA Conjecture} holds for $n \leq 9$.
Although the conjecture remains unresolved,
one of its two directions has already been established:

\begin{rem}\label{Easy direction}
{\rm \cite[p. 43]{BBH}}
\rm
In Conjecture \ref{NAA Conjecture}, (b) implies (a).
\end{rem}

Remark \ref{Easy direction} implies that
the family of nonsingular symmetric 0--1 matrices whose
inverses are equimodular with constant diagonal is a subset of
the family of 
nonsingular symmetric 0--1 matrices that end with $\tt NAA$.
Conjecture \ref{NAA Conjecture} speculates that
these two families are one and the same.
Furthermore, Remark \ref{Easy direction} implies that,
to establish Conjecture \ref{NAA Conjecture} (if it is true),
it suffices to show that (a) implies (b);
i.e., it suffices to show that the inverse of an arbitrary
symmetric 0--1 matrix that ends with $\tt NAA$ is
equimodular with constant diagonal;
to show this,
it suffices to show that
all of the principal minors of order $n-1$ of
an arbitrary symmetric matrix in $\zonn$ that
ends with $\tt NAA$ are equal:

\begin{rem}\label{(n-1)-minors}
{\rm \cite[p. 45]{BBH}}
\rm
Let $A \in \zonn$ be a symmetric matrix that ends with $\tt NAA$.
If all of the principal minors of $A$ of order $n-1$ are equal, then
$A^{-1}$ is equimodular with constant diagonal.
\end{rem}


One of our chief objectives is to establish that,
if the matrix $A$ in Conjecture \ref{NAA Conjecture} satisfies
one particular condition, then the conjecture holds;
we provide more details regarding 
this condition at the end of this section,
immediately after additional terminology is introduced.
We shall also cast more light upon the
two families of nonsingular symmetric 0--1 matrices that 
have been conjectured 
(in Conjecture \ref{NAA Conjecture}) to be one and the same:
The family whose members end with $\tt NAA$ and
the family whose members are nonsingular matrices whose
inverses are equimodular with constant diagonal.

For us, all graphs are simple, undirected and loopless.
The cycle graph on $n \geq 3$ vertices is denoted by $C_n$, and
$C_n$ is \textit{odd} (respectively, \textit{even}),
if $n$ is odd (respectively, even).
For a given graph $G$,
its number of vertices is its \textit{order};
the \textit{odd-girth} of $G$ is defined to be the length
(number of vertices, or edges) of a shortest odd cycle in $G$,
and, if no such cycle exists, then it is defined to be $\infty$;
and the 0--1 adjacency matrix of $G$ is denoted by $A(G)$.
A given square matrix all of whose diagonal entries are zero is
said to have ``zero diagonal.''
With a given symmetric matrix
$B = [b_{ij}] \in \Rnn$ that has zero diagonal,
we associate a graph, which we denote by $G(B)$,
whose vertex set is $[n]$, and
whose edges are determined as follows:
The graph $G(B)$ has an edge between the vertices $i,j \in [n]$
if and only if
$b_{ij}$ is nonzero.

For a given matrix $B \in \Rnn$ with a
nonsingular principal submatrix $B[\mu]$, recall that
the Schur complement of $B[\mu]$ in $B$ is the matrix
$B/B[\mu] :=
B[\mu^c] - B[\mu^c,\mu](B[\mu])^{-1}B[\mu,\mu^c]$,
where $\mu^c = [n] \setminus \mu$ (see, for example, \cite{Schur}).
The Kronecker product  (also known as tensor product) of
two given matrices $B$ and $C$ is denoted by $B \otimes C$.
The matrices $B$ and $C$ are said to be
{\em permutationally similar} if there exists a
permutation matrix $P$ such that $C=P^TBP$.
By $I_n$, $O_{m,n}$ and $J_n$ we denote, respectively,
the identity matrix of order $n$,
the $m \times n$ zero matrix and
the $n \times n$ matrix all of whose entries are equal to $1$;
furthermore, $O_{n,n}:=O_n$.

All of the matrices treated in this paper are real
(i.e., they are over $\R$).
However,
we shall, at one point in Section \ref{s:NAA},
analyze a given (real) 0--1 matrix by treating it as a matrix over
the field of order $2$, which we shall denote by $\ZZ$.

The purpose of Section \ref{s:known epr} is
merely to list some known results that
will be cited in the subsequent pages.
In Section \ref{s:AN and N},
our focus is on epr-sequences that start with $\tt AN$
(which are of particular relevance to
Conjecture \ref{NAA Conjecture}),
as well as on epr-sequences that start with $\tt N$:
The epr-sequences that
start with {\tt AN} but whose fourth term is not $\tt S$ are
completely characterized;
and it is shown that,
if $B$ is an arbitrary symmetric matrix with zero diagonal
(i.e., a matrix that starts with $\tt N$),
then the odd-girth of $G(B)$,
as well as whether or not $G(B)$ is bipartite or an odd cycle,
is deducible from $\epr(B)$.
Section \ref{s:NAA} is devoted to
symmetric 0--1 matrices that end with $\tt NAA$
(i.e., to the matrices in the first of the two
families associated with Conjecture \ref{NAA Conjecture});
in particular, it is established that
Conjecture \ref{NAA Conjecture} holds if $A$ does not have
both a zero and a nonzero principal minor of order $n-4$
(if $n \geq 5$);
moreover, motivated by Remark \ref{(n-1)-minors},
attention is paid to the parity of
the principal minors of these matrices, especially to
the parity of their principal minors of order $n-1$ and
the parity of their determinants,
with the latter of which shown to all be even.
Section \ref{s:equim} focuses on 
(not necessarily symmetric)
nonsingular 0--1 matrices whose inverses are equimodular,
due to their relation with the second of the two
families associated with Conjecture \ref{NAA Conjecture}
(namely, 
the family of nonsingular symmetric 0--1 matrices whose
inverses are equimodular with constant diagonal);
for an arbitrary (not necessarily symmetric)
nonsingular matrix $A \in \zonn$ with $n\geq 3$,
we establish necessary conditions for 
$A^{-1}$ to be equimodular; 
examples of such conditions are the following:
each row and column of $A$ has an 
even number of nonzero entries;
each entry of $A^{-1}$ is the reciprocal of an even integer;
$\det(A)$ is even;
the difference between any two rows of $A$, 
as well as the difference between any two columns of $A$, 
has an even number of nonzero entries;
if $A$ is symmetric, then $A$ has an 
even number of nonzero diagonal entries;
if $A$ is symmetric and 
$\vec{a}_k$ is the $k$th column of $A$, 
then $A-\vec{a}_k\vec{a}_k^T$
has an even number of nonzero diagonal entries.

\section{Some known results}\label{s:known epr}
$\null$
\indent
The purpose of this section is to list some known results
(of which most are about epr-sequences) that will be
cited in the remaining pages;
some of them are assigned names by which they
will be referenced.
Results that are restricted to
the field $\R$ may hold over other fields
(most results are stated in
the context of $\R$ merely for simplicity).
The first result listed is well-known (and elementary):

\begin{thm}\label{adjoint-result}
Let $B \in \Rnn$ be nonsingular and $B^{-1} = [\beta_{ij}]$.
Then, for all $i,j \in [n]$,
\[
\beta_{ij} = (-1)^{i+j}\frac{\det(B(\{j\},\{i\}))}{\det(B)}
\quad \text{ and } \quad
\beta_{ii} = \frac{\det(B(\{i\}))}{\det(B)}.
\]
\end{thm}

The next result, which is also well-known, 
is concerned with Schur complements.

\begin{thm}\label{Schur}
{\rm  (Schur Complement Theorem.)}
Let $B \in \Rnn$,
$\rank(B)=r$,
$\mu \subset [n]$, 
$\mu^c = [n]\setminus \mu$ and $|\mu| = k$.
Suppose that $B[\mu]$ is nonsingular,
and let $C = B/B[\mu]$.
Then the following statements hold:
\begin{enumerate}
\item [$(i)$] 
$C$ is an $(n-k)\times (n-k)$ matrix.
\item[$(ii)$]
If $B$ is symmetric, then $C$ is symmetric.
\item [$(iii)$] {\rm \cite[p. 14]{Schur}}
$\rank(C) = r-k$.
\item[$(iv)$] {\rm \cite[Theorem 1.2]{Schur}}
If $B$ is nonsingular, then
$C$ is nonsingular and $C^{-1} = B^{-1}[\mu^c]$.
\item [$(v)$] {\rm \cite[p. 771]{Brualdi & Schneider}}
If the indexing of $C$ is inherited from $B$, then,
for all $\gamma \subseteq \mu^c$,
\[
\det(C[\gamma]) = 
\frac{\det(B[\gamma \cup \mu])}{\det(B[\mu])}.
\]
\end{enumerate}
\end{thm}

The following is also a well-known fact
(see, for example, \cite[Theorem 1.1]{BIRS13}),
and it states that the rank of an 
arbitrary symmetric matrix $B$ is equal to
the order of a largest nonsingular principal submatrix of $B$,
which is why, in Section \ref{s:intro}, we adopted the 
convention of calling the rank of a symmetric matrix ``principal.''

\begin{thm}\label{thm: rank of a symm mtx}
Let $B \in \Rnn$ be symmetric.
Then \ $\rank(B) = \max \{ |\mu| : \det (B[\mu]) \neq 0 \}$,
where the maximum over the empty set is defined to be 0.
\end{thm}


The next result follows from
the Schur Complement Theorem.

\begin{thm}\label{SchurAN}
{\rm \cite[Corollary 1.11]{EPR-Hermitian}}
Let $B \in \Rnn$ be symmetric,
$\epr(B)=\ell_1 \ell_2 \cdots \ell_n$,
$\mu \subset [n]$ and $|\mu| = k$.
Suppose that $B[\mu]$ is nonsingular, and let
$C = B/B[\mu]$ and
$\epr(C)=\ell'_{1} \ell'_2 \cdots \ell'_{n-k}$.
If $\ell_{j+k} \in \{{\tt A,N}\}$,
then $\ell'_j=\ell_{j+k}$,  for all $j\in [n-k]$.
\end{thm}

The epr-sequence of the inverse of an arbitrary nonsingular
symmetric matrix is deducible from that of the given matrix:

\begin{thm}\label{Inverse Thm}
{\rm\cite[Theorem 2.4]{EPR} (Inverse Theorem.)}
Let $B \in \Rnn$ be symmetric.
Suppose that $B$ is nonsingular.
If $\epr(B) = \ell_1\ell_2 \cdots \ell_{n-1}\tt{A}$, then
$\epr(B^{-1}) = \ell_{n-1}\ell_{n-2} \cdots \ell_{1}\tt{A}$.
\end{thm}

The next result allows one to 
deduce some of the terms in the epr-sequence of 
a principal submatrix of an arbitrary symmetric matrix, 
provided that the epr-sequence of the latter is known.

\begin{thm}\label{Inheritance}
{\rm{\cite[Theorem 2.6]{EPR}}}
{\rm (Inheritance Theorem.)}
Let $B \in \Rnn$ be symmetric, $m \leq n$ and $1\le i \le m$.
Then the following statements hold:
\ben
\item
If $[\epr(B)]_i={\tt N}$, then  $[\epr(C)]_i={\tt N}$,
for all $m\times m$ principal submatrices $C$.
\item If  $[\epr(B)]_i={\tt A}$, then  $[\epr(C)]_i={\tt A}$,
for all $m\times m$ principal submatrices $C$.

\item If $[\epr(B)]_m={\tt S}$, then there exist $m\times m$ principal submatrices $C_A$ and $C_N$ of $B$ such that $[\epr(C_A)]_m = {\tt A}$ and $[\epr(C_N)]_m = {\tt N}$.
\item If $i < m$ and $[\epr(B)]_i = {\tt S}$, then there exists an $m \times m$ principal submatrix $C_S$ such that $[\epr(C_S)]_i ={\tt S}$.
\een
\end{thm}

If the epr-sequence of an arbitrary 
$n \times n$ symmetric matrix $B$ is known,
then there is a way to append a row and a column to $B$ so that
the epr-sequence of the resulting $(n+1) \times (n+1)$ matrix ``inherits'' the first letter of $\epr(B)$ as well as any $\tt N$s (and $\tt S$s) in $\epr(B)$:

\begin{obs}\label{Obs: Append Ns and preserve ell1}
{\rm\cite[Observation 2.19]{EPR}}
Let $B = [b_{ij}] \in \Rnn$ be symmetric and
$\epr(B) = \ell_1\ell_2\cdots\ell_n$.
Let $\vec{x}$ be the $n$th column of $B$,
\[
B' =
\left(
\begin{array}{c|c}
B & \vec{x} \\
\hline
\vec{x}^T & b_{nn}
\end{array}
\right)
\]
and $\epr(B') = \ell'_1\ell'_2 \cdots \ell'_{n+1}$.
Then $\epr(B')=\ell_1\ell_2'\cdots\ell_n'{\tt N}$.
Moreover, for $2\le j\le n$, the following statements hold:
If $\ell_j={\tt N}$, then $\ell_j'={\tt N}$; and,
if $\ell_j \in \{\tt A,S\}$, then $\ell_j'={\tt S}$.

\end{obs}

If there are two consecutive $\tt N$s in
the epr-sequence of a symmetric matrix,
then each letter in the sequence from that point forward is $\tt N$:

\begin{thm}\label{NN result}
{\rm \cite[Theorem 2.3]{EPR} ($\tt NN$ Theorem.)}
Let $B \in \Rnn$ be symmetric and
$\epr(B) = \ell_1\ell_2 \cdots \ell_n$.
Suppose that $\ell_k = \ell_{k+1} = \tt{N}$, for some $k$.
Then, for all $j \geq k$, $\ell_j = \tt{N}$.
\end{thm}


\begin{thm}\label{NSA-NAN-NAS}
{\rm \cite[Corollary 2.7 \& Theorem 2.14]{EPR}}
Let $B \in \Rnn$ be symmetric.
Then none of the following are subsequences of $\epr(B)$:
\[
{\tt NSA}; \quad
{\tt NAN}; \quad
{\tt NAS}.
\]
\end{thm}

We shall reference Theorem \ref{NSA-NAN-NAS} by saying,
for example, that ``{\tt NSA} is forbidden''
(likewise with {\tt NAN} and {\tt NAS}).

We conclude this section by listing some
other known results about epr-sequences.

\begin{thm}\label{ASN...A...}
{\rm\cite[Corollary 2.7]{EPR}}
Let $B \in \Rnn$ be symmetric and 
$\epr(B) = \ell_1\ell_2 \cdots \ell_n$.
Suppose that, for some $k \in [n-2]$,
$\ell_k\ell_{k+1}\ell_{k+2} = \tt ASN$.
If $j \geq k+3$,  then $\ell_j \neq \tt A$.
\end{thm}



\begin{thm}\label{ANA}
{\rm \cite[Theorem 2.6]{XMR-Classif}}
Let $B \in \Rnn$ be symmetric and
$\epr(B) = \ell_1\ell_2 \cdots \ell_n$.
Suppose that
$\ell_1\ell_2 \cdots \ell_{n-1}$ contains $\tt ANA$ as a subsequence.
Then $\epr(B)$ is of the form $\tt{\overline{A}ANAA\overline{A}}$.
\end{thm}

\begin{thm}\label{SNA}
{\rm \cite[Theorem 2.7]{XMR-Classif}}
Let $B \in \Rnn$ be symmetric and
$\epr(B) = \ell_1\ell_2 \cdots \ell_n$.
Then $\tt{SNA}$ is not a subsequence of
$\ell_1\ell_2 \cdots \ell_{n-2}$.
\end{thm}

\begin{thm}\label{ANSNSN(SN)(N)}
{\rm \cite[Proposition 3.5]{XMR-Classif}}
Let $B \in \Rnn$ be symmetric with
$\epr(B) = {\tt ANSNSN}\ell_7 \cdots \ell_n$.
Then, for $k \geq 1$, $\ell_{2k} = \tt{N}$.
Furthermore, $\epr(B)$ is of one of the following forms:
$\tt{ANSNSN\overline{SN}}\hspace{0.04cm}\overline{\tt{N}}$
or
${\tt ANSNSN\overline{SN}A}$.
\end{thm}


\begin{obs}\label{epr of graphs}
{\rm\cite[Observation 3.1]{EPR}}
If $n \geq 3$ is odd, then $\epr(A(C_n))=\tt \OL{NS}NAA$.
\end{obs}

Observe that the matrix $A$ in (\ref{5x5-matrix-example}) is $A(C_5)$.

\section{Epr-sequences that start with {\tt AN} or {\tt N}}\label{s:AN and N}
$\null$
\indent
In this section,
results of particular relevance to our results concerning
Conjecture \ref{NAA Conjecture} are established.
We start by focusing on epr-sequences that start with $\tt AN$;
the relevance of these is made evident by
the following observation
(the first statement of which follows from the Inverse Theorem, while second statement follows readily).

\begin{obs}
Let $n \geq 2$ be an integer and $B, C \in \Rnn$  be symmetric.
If $B$ ends with $\tt NAA$, then $B^{-1}$ starts with $\tt AN$.
If $C \neq O_n$ and is equimodular with constant diagonal,
then $C$ starts with $\tt AN$.
\end{obs}

We need the next lemma, which generalizes \cite[Proposition 2.13(4)]{EPR}.

\begin{lem}\label{lemma from Schur}
Let $B \in \Rnn$ be symmetric, $\rank(B)=r$,
$\mu \subseteq [n]$ and $|\mu| = k$.
Suppose that $B[\mu]$ is nonsingular.
If $m$ is an integer with $k< m \leq r$,
then $B[\mu]$ is contained in either an
$m \times m$ or an  $(m+1) \times (m+1)$
nonsingular principal submatrix of $B$.
\end{lem}

\bpf
Suppose that $m$ is an integer with $k< m \leq r$.
Let $C = B/B[\mu]$, and
assume that the indexing of $C$ is inherited from $B$.
By the Schur Complement Theorem,
$C \in \R^{(n-k) \times (n-k)}$, $C$ is symmetric and
$\rank(C) = r-k$.
Let $\epr(C) = \ell_1\ell_2 \cdots \ell_{n-k}$.
Then, as the rank of $C$ is principal, the $\tt NN$ Theorem implies that
$\ell_{m-k} \neq \tt N$ or $\ell_{m+1-k} \neq \tt N$
(otherwise, $\rank(C) < m-k \leq r-k$, which is a contradiction).
Thus, for some $\gamma \subseteq [n]\setminus \mu$ with
$|\gamma| = m-k$ or $|\gamma| = m+1-k$,
$\det(C[\gamma]) \neq 0$.
By the Schur Complement Theorem,
\[
\det(C[\gamma]) =
\frac{\det(B[\gamma \cup \mu])}{\det(B[\mu])},
\]
implying that $\det(B[\gamma \cup \mu]) \neq 0$.
Then, as $|\gamma \cup \mu| \in \{m,m+1\}$,
the desired conclusion follows.
\epf

A complete characterization of the epr-sequences that
start with $\tt AN$ but whose fourth letter is not $\tt S$ may
now be established.

\begin{thm}\label{ANSN}
Let $\sigma = {\tt AN}\ell_3\ell_4 \cdots \ell_n$ be a
sequence from $\{\tt A,N,S\}$ with $\ell_4 \neq \tt S$.
Then there exists a symmetric matrix $B \in \Rnn$ with
$\epr(B) = \sigma$ if and only if $\sigma$ is of one of the following forms:
\[
{
\tt AN\OL{N}, \quad
\tt ANA\OL{A}, \quad
\tt ANSN\OL{SN}\hspace{0.04cm}\OL{N}, \quad
\tt ANSN\OL{SN}A, \quad
\tt ANSNAA, \text{ or } \
\tt ANSNSSN\OL{N}.
}
\]
\end{thm}

\bpf
Suppose that
there exists a symmetric matrix $B \in \Rnn$ with $\epr(B) = \sigma$.
As it is straightforward  to verify that  the desired conclusion holds for $n \leq 3$
(by referring to \cite[Table 1]{EPR}, for example),
we assume that $n \geq 4$.
If $\ell_3 = \tt N$, then the $\tt NN$ Theorem implies that
$\sigma = \tt ANNN\OL{N}$.
If $\ell_3 = \tt A$, then Theorem \ref{ANA} implies that
$\sigma = \tt ANAA\OL{A}$.
Assume that $\ell_3 = \tt S$.
By the hypothesis, $\ell_4 \in \{\tt A,N\}$.
Then, as $\tt NSA$ is forbidden, $\ell_4 = \tt N$,
implying that $\sigma$ starts with $\tt ANSN$.
Since it is straightforward  to verify that
the desired conclusion holds for $n \leq 7$
(by referring to \cite[Table 1]{EPR}, for example),
we assume that $n \geq 8$.
By Theorem \ref{SNA},  $\ell_5 \neq \tt A$.
We consider two cases, based on $\ell_6$.

\noindent
\textbf{Case 1}: $\ell_6 = \tt N$.

\noindent
If $\ell_5 = \tt N$, then the $\tt NN$ Theorem implies that
$\epr(B)$ is of the form 
$\tt ANSN\OL{SN}\hspace{0.04cm}\OL{N}$
(with $\OL{\tt SN}$ vacuous).
If $\ell_5 = \tt S$, then it follows from
Theorem \ref{ANSNSN(SN)(N)} that
$\epr(B)$ is of one of the following two forms:
$\tt ANSNSN\OL{SN}\hspace{0.04cm}\OL{N}$ or
$\tt ANSNSN\OL{SN}A$.

\noindent
\textbf{Case 2}: $\ell_6 \neq \tt N$.

\noindent
It follows from the $\tt NN$ Theorem that $\ell_5 \neq \tt N$.
Then, as $\ell_5 \neq \tt A$, $\ell _5 = \tt S$.
Thus, $\ell_5\ell_6 \in \{\tt SA, SS\}$.
As $\tt NSA$ is forbidden, $\ell_5\ell_6 = \tt SS$.
Thus,
$\epr(B) = {\tt ANSNSS}\ell_7 \cdots \ell_n$.
Hence, $\rank(B) \geq 6$.
We now show that $\rank(B)=6$.
Suppose on the contrary that $\rank(B) \geq 7$.
Let $B[\mu]$ be a $6 \times 6$ nonsingular principal submatrix of $B$.
By Lemma \ref{lemma from Schur},
$B[\mu]$ is contained in either a
$7 \times 7$ or an 
$8 \times 8$ nonsingular principal submatrix of $B$.
Thus, there exists $\omega \subseteq [n]\setminus \mu$ with
$1 \leq |\omega| \leq 2$ such that $B[\mu \cup \omega]$ is nonsingular.
Let $\epr(B[\mu \cup \omega]) = \ell'_1\ell'_2 \cdots \ell'_{6 + |\omega|}$.
By the Inheritance Theorem,
$\ell'_1\ell'_2\ell'_3\ell'_4 = {\tt AN}\ell'_3{\tt N}$.
The nonsingularity of $B[\mu \cup \omega]$ implies that
$\ell'_{6 + |\omega|} = \tt A$.
Thus, the {\tt NN} Theorem implies that $\ell'_3 \neq \tt N$.
Then, as {\tt NAN} is forbidden, $\ell'_3 = \tt S$.
Thus, either
$\epr(B[\mu \cup \omega]) = {\tt ANSN}\ell'_5 \ell'_{6}{\tt A}$
or
$\epr(B[\mu \cup \omega]) = {\tt ANSN}\ell'_5\ell'_{6}\ell'_7 {\tt A}$.
By \cite[Table 1]{EPR},
$\epr(B[\mu \cup \omega]) \neq {\tt ANSN}\ell'_5\ell'_{6}\ell'_7 {\tt A}$.
Hence,
$\epr(B[\mu \cup \omega]) = {\tt ANSN}\ell'_5 \ell'_{6}{\tt A}$.
It follows from  \cite[Table 1]{EPR} that
$\epr(B[\mu \cup \omega]) = {\tt ANSNSNA}$.
Thus,
each $6 \times 6$ principal submatrix of $B[\mu \cup \omega]$ is singular.
Then, as the $6 \times 6$ matrix $B[\mu]$ is a principal submatrix of $B[\mu \cup \omega]$, $B[\mu]$ is singular, which is a contradiction.
We conclude that $\rank(B) = 6$.
Then, as the rank of $B$ is principal,
$\ell_j = \tt N$ for all $j \geq 7$,
implying that $\epr(B)$ is of the form $\tt ANSNSSNN\OL{N}$.

We now establish the other direction.
If $\sigma = \tt AN\OL{N}$,
then $\epr(J_n) = \sigma$ (that follows readily).
If $\sigma = \tt ANA\OL{A}$,
then $\epr(J_n-2I_n)=\sigma$
(see \cite[Proposition 2.5]{XMR-Classif}).
If $\sigma = \tt ANSNAA$,
then \cite[Table 1]{EPR} implies that there exists a
symmetric matrix $C^{(0)} \in \R^{6 \times 6}$ with
$\epr(C^{(0)}) = \sigma$.
Assume that $\sigma$ is of the form $\tt ANSNSSN\OL{N}$.
By applying 
Observation \ref{Obs: Append Ns and preserve ell1} to 
$C^{(0)}$, we conclude that there exists a symmetric matrix
$C^{(1)} \in \R^{7 \times 7}$ with 
$\epr(C^{(1)}) = \tt ANSNSSN$.
Similarly, by applying 
Observation \ref{Obs: Append Ns and preserve ell1} to 
$C^{(1)}$, we conclude that there exists a symmetric matrix $C^{(2)} \in \R^{8 \times 8}$ with
$\epr(C^{(2)}) = \tt ANSNSSNN$.
By continuing in this manner, one concludes that 
there exists a symmetric matrix $C^{(n-6)} \in \Rnn$ with 
$\epr(C^{(n-6)}) = \tt ANSNSSN\OL{N} =\sigma$
(with $\tt \OL{N}$ having $n-6$ copies of $\tt N$).
Assume that $\sigma$ is of the form $\tt ANSN\OL{SN}A$.
By Observation \ref{epr of graphs} and the Inverse Theorem,
$\epr((A(C_n))^{-1})= \sigma$.
Finally, assume that $\sigma$ is of the form
$\tt ANSN\OL{SN}\hspace{0.04cm}\OL{N}$.
Because of Observation \ref{Obs: Append Ns and preserve ell1},
it suffices to establish the desired conclusion in the case where
$\tt \OL{N}$ is assumed to be vacuous.
If $\tt \OL{SN}$ is vacuous
(meaning that $\sigma = \tt ANSN$), then,
by \cite[Table 1]{EPR}, there exists a symmetric matrix
$M \in \Rnn$ with $\epr(M) = \tt \sigma$.
Assume that $\tt \OL{SN}$ is not vacuous
(meaning that $\sigma = \tt ANSNSN\OL{SN}$).
Let $N \in \R^{(n-1)\times (n-1)}$ be a symmetric matrix with
$\epr(N) = \tt ANSN\OL{SN}A$
(we already showed that $N$ exists).
Then, by Observation \ref{Obs: Append Ns and preserve ell1},
there exists a symmetric matrix $N' \in \Rnn$ with
$\epr(N') = \sigma$.
\epf

For the remainder of this section,
our attention is turned to epr-sequences that start with $\tt N$.
A natural question is the following one:
If $B$ is an arbitrary symmetric matrix with zero diagonal
(i.e., a matrix that starts with $\tt N$),
what properties of $G(B)$ are deducible from $\epr(B)$?
Its odd-girth is one such property:

\begin{thm}\label{odd-girth result}
Let $B \in \Rnn$ be symmetric and
$\epr(B) = \ell_1\ell_2 \cdots \ell_n$.
Suppose that $\ell_1 = \tt N$ and
let $g$ be the odd-girth of \ $G(B)$.
Then, for all odd $j$ in the interval $[1,n] \cap [1,g)$, $\ell_j = \tt N$.
Moreover,
$g=\min \{j \in[n] : j \text{ is odd and } \ell_j \neq \tt N \}$,
where the minimum over the empty set is defined to be $\infty$.
In particular, if $G(B)$ is neither bipartite nor an isolated vertex,
then $g$ is equal to the order of a
smallest nonsingular principal submatrix of $B$ whose order is odd.
\end{thm}

\bpf
If $n=1$, then $G(B)$ is an isolated vertex,  and, therefore,
$g=\infty$ and all of the desired conclusions follow readily.
Now assume that $n \geq 2$.
Let $t = \min \{j \in[n] : j \text{ is odd and } \ell_j \neq \tt N \}$.
We divide the proof into two cases:

\noindent
\textbf{Case 1}: $G(B)$ is bipartite.

\noindent
It follows, then, that $g = \infty$.
Thus, $[1,n] \cap [1,g) = [1,n]$.
By the hypothesis, $\ell_1 = \tt N$.
Since each principal submatrix of $B$ whose order is at least two is bipartite,
each principal submatrix of $B$ whose order is odd is singular,
implying that $\ell_j = \tt N$, for all odd $j \in [1,n]$, as desired.
Hence, $t = \infty =g$, as desired.

\noindent
\textbf{Case 2}: $G(B)$ is not bipartite.

\noindent
Since $G(B)$ is not an isolated vertex, it has an odd cycle,
implying that $g \geq 3$.
Thus, $[1,n] \cap [1,g) = [1, g-1]$.
Let $j \in [1, g-1]$ be odd.
If $j=1$, then, by the hypothesis, $\ell_j = \tt N$.
Assume that $j \geq 3$.
Then, as $j<g$, the graph of each $j \times j$ principal submatrix of $B$
does not have any odd cycles, implying that such graphs are all
bipartite and of odd order.
Thus, $\ell_j = \tt N$.
It follows that,
for all odd $j$ in the interval $[1, g-1]$, $\ell_j = \tt N$, as desired.

We now show (still within Case 2) that $t=g$.
It follows from the final conclusion of 
the previous paragraph that $t \geq g$.
Suppose on the contrary that $t>g$.
Then, as $g$ is odd, $\ell_g = \tt N$.
Let $B[\mu]$ be a $g \times g$ (principal) submatrix whose graph,
$G(B[\mu])$, has a cycle of length $g$
(such submatrix exists because $g$ is the odd-girth of $G(B)$).
Since every cycle of $G(B[\mu])$ is a cycle of $G(B)$,
and because $g$ is the odd-girth of $G(B)$,
$G(B[\mu])$ does not have any odd cycles whose
length is less than $g$.
Then, as $g$ is odd and $G(B[\mu])$ is of order $g$,
$G(B[\mu])$ does not have any other edges
besides those that are part of its cycle of length $g$
(otherwise, $G(B[\mu])$ would have a shorter odd cycle).
Hence, $G(B[\mu]) = C_{g}$.
It follows (from well-known facts) that $B[\mu]$ is nonsingular,
implying that $\ell_{g} \neq \tt N$, which is a contradiction.
We conclude that $t=g$, as desired.

Since a graph whose odd-girth is $\infty$ is
either bipartite or an isolated vertex,
the last statement of the theorem is readily deducible from
the one preceding it.
\epf

As a consequence of
Theorem \ref{odd-girth result}, the following statement holds:
If $G$ is a non-bipartite graph of order $n \geq 2$ with odd-girth $g$,
then $g$ is equal to the order of a
smallest nonsingular principal submatrix of $A(G)$ whose order is odd.

\begin{lem}\label{lemma for bipartite and cycle}
Let  $B \in \Rnn$ be symmetric,
$\epr(B) = \ell_1\ell_2 \cdots \ell_n$ and
$k$ be an odd integer with $1 \leq k \leq n$.
Suppose that, for all odd integers $j$ with $1 \leq j \leq k$, $\ell_j = \tt N$.
If $1 \leq k <\rank(B)$,
then $\ell_1\ell_2 \cdots \ell_k$ is of the form $\tt N\OL{SN}$.
\end{lem}

\bpf
Suppose that $1 \leq k < \rank(B)$.
If $k = 1$, then
the desired conclusion follows by the hypothesis.
Assume that $k \geq 3$.
It suffices to show that,
for all even integers $i$ with $2 \leq i \leq k-1$, $\ell_i = \tt S$.
Let $i$ be an even integer with $2 \leq i \leq k-1$.
As $\rank(B) > k$, and because the rank of $B$ is principal,
the $\tt NN$ Theorem implies that $\ell_i \neq \tt N$
(otherwise, $\rank(B) < i < k$).
Then, as $\tt NAN$ is prohibited, $\ell_i = \tt S$.
\epf

The following fact generalizes
some of the results in \cite[Observation 3.1]{EPR},
and it implies that the following statement holds:
If $B$ is an arbitrary symmetric matrix with zero diagonal
(i.e., a matrix that starts with $\tt N$),
then whether or not $G(B)$ is
bipartite or an odd cycle is deducible from $\epr(B)$.

\begin{thm}\label{bipartite-cycle result}
Let $B \in \Rnn$ be symmetric and $\epr(B) = \ell_1\ell_2 \cdots \ell_n$.
Suppose that $\ell_1 = \tt N$.
Then the following statements hold:
\ben
\item $G(B)$ is bipartite
if and only if
$\epr(B)$ is of one of the following forms:
$\tt NN\OL{N}$ or
$\tt \OL{NS}NA$ or
$\tt NS\OL{NS}N\OL{N}$.
\item $G(B) = C_{2k+1}$ for some integer $k \geq 1$
if and only if
$\epr(B) = \tt \OL{NS}NAA$.
\een
\end{thm}

\bpf
We start by establishing statement (1).
Suppose that $G(B)$ is bipartite.
Thus, $n \geq 2$.
If $B = O_n$, then, clearly,
$\epr(B)$ is of the form $\tt NN\OL{N}$.
Now assume that $B \neq O_n$.
Then, as $B$ is symmetric and $\ell_1 = \tt N$, $\rank(B) \geq 2$.
We now show that $\epr(B)$ is of one of the following forms:
$\tt \OL{NS}NA$ or $\tt NS\OL{NS}N\OL{N}$.
As $G(B)$ is bipartite,
the odd-girth of $G(B)$ is equal to $\infty$.
By Theorem \ref{odd-girth result},
for all odd $j$ in the interval $[1,n] \cap [1,\infty)=[1,n]$, $\ell_j = \tt N$.
Thus, $\rank(B)$ is even (and nonzero).
Let $k=\rank(B)-1$.
Observe that $k$ is odd and 
$1 \leq k \leq n-1$ and $1 \leq k < \rank(B)$.
By Lemma \ref{lemma for bipartite and cycle},
$\ell_1\ell_2 \cdots \ell_k$ is of the form $\tt N\OL{SN}$,
which is of the form $\tt \OL{NS}N$.
As the rank of $B$ is principal, and because $k+1=\rank(B)$,
$\ell_{k+1} \in \tt \{A,S\}$ and 
$\ell_q = \tt N$ for all $q \geq k+2$.
Observe that
$\ell_{k+1}$ is the only term of $\epr(B)$ that remains unknown.
If $k=n-1$, then
$\epr(B) = \ell_1\ell_2 \cdots \ell_{k+1}$
and, therefore, $\ell_{k+1}= \tt A$
(because the epr-sequence of no matrix ends with $\tt S$),
implying that $\epr(B)$ is of the form $\tt \OL{NS}NA$.
Now assume that $k \leq n-2$.
Then, as $\tt NAN$ is prohibited,
$\ell_{k+1} = \tt S$ and, therefore,
$\epr(B)$ is of the form
$\tt \OL{NS}NSN\OL{N}$, which is of the form
$\tt NS\OL{NS}N\OL{N}$.

For the other direction,
suppose that $\epr(B)$ is of one of the following forms:
$\tt NN\OL{N}$ or $\tt \OL{NS}NA$ or $\tt NS\OL{NS}N\OL{N}$.
Thus, $n \geq 2$.
Moreover, by Theorem \ref{odd-girth result},
the odd-girth of $G(B)$ is $\infty$.
Then, as $G(B)$ is not an isolated vertex (because $n \geq 2$),
$G(B)$ is bipartite.

We now establish statement (2).
Suppose that $G(B) = C_{2k+1}$ for some integer $k \geq 1$.
Thus, $n=2k+1$ (which is an odd integer) and $n \geq 3$.
Observe that the odd-girth of $G(B)$ is $n$.
Then, by Theorem \ref{odd-girth result}, $\ell_n \neq \tt N$.
Then, as the epr-sequence of no matrix ends with $\tt S$,
$\ell_n = \tt A$.
By Theorem \ref{odd-girth result},
for all odd integers $j$ in the interval
$[1,n] \cap [1,n)=[1,n)$, $\ell_j = \tt N$.
Since $n-2$ is odd, $\ell_{n-2} = \tt N$.
Then, as $\ell_n =\tt A$,
the $\tt NN$ Theorem implies that $\ell_{n-1} \neq \tt N$.
It follows from the fact that $\tt NSA$ is prohibited
that $\ell_{n-1} = \tt A$.
It remains to show that $\ell_i = \tt S$ for
all even integers with $2 \leq i \leq n-3$;
that follows by applying 
Lemma \ref{lemma for bipartite and cycle} with $k=n-2$.

For the other direction,
suppose that $\epr(B) = \tt \OL{NS}NAA$.
Thus, $n \geq 3$ and $n$ is odd.
Hence, it suffices to show that $G(B) = C_n$.
Let $g$ be the odd-girth of $G(B)$.
By Theorem \ref{odd-girth result}, $g=n$.
Thus, $C_{n}$ is a subgraph of $G(B)$.
Then, as $G(B)$ is of order $n$, and because $n$ is odd,
$G(B)$ does not have any other edges besides those that
are part of its cycle of length $n$
(otherwise, $G(B)$ would have an
odd cycle that is shorter than $C_n = C_g$).
Hence, $G(B) = C_n$, as desired.
\epf

\section{Symmetric 0--1 matrices that end with {\tt NAA}}\label{s:NAA}
$\null$
\indent
This section is devoted to
symmetric 0--1 matrices that end with $\tt NAA$
(i.e., to the matrices in the first of the two
families associated with Conjecture \ref{NAA Conjecture}).

The particular relevance of the results of the previous section is
elucidated by the proof of the next theorem.

\begin{thm}\label{NSNAA-SSNAA}
Let $n \geq 5$ and $A \in \zonn$ be symmetric.
Suppose that $A$ ends with $\tt NAA$.
Then $A$ ends with either $\tt NSNAA$ or $\tt SSNAA$.
Moreover, if $A$ ends with $\tt NSNAA$, then
$n$ is odd and $A$ is permutationally similar to $A(C_n)$.
Consequently, if $n$ is even, then $A$ ends with $\tt SSNAA$.
Furthermore, if $A$ is not permutationally similar to $A(C_n)$,
then $A$ ends with $\tt SSNAA$.
\end{thm}

\bpf
Let $\epr(A) = \ell_1\ell_2 \cdots \ell_n$.
Thus, $\ell_{n-2}\ell_{n-1}\ell_n = \tt NAA$.
By the $\tt NN$ Theorem, $\ell_{n-3} \neq \tt N$.
We claim that $\ell_{n-3} = \tt S$.
Suppose on the contrary that $\ell_{n-3} = \tt A$.
By Theorem \ref{ANA}, $\epr(A) = \tt \OL{A}AANAA$.
Thus, $A$ is a 0--1 matrix that starts with $\tt AA$,
implying that $A = I_n$ and, therefore, that
$\epr(A) = \tt AAAAA\OL{A}$, which is a contradiction.
It follows that our claim is true; that is, $\ell_{n-3} = \tt S$.
By Theorem \ref{ASN...A...}, $\ell_{n-4} \neq \tt A$.
Thus, $A$ ends with either $\tt NSNAA$ or $\tt SSNAA$.

Suppose that $A$ ends with $\tt NSNAA$.
Then, as $A$ is a nonsingular 0--1 matrix,
$A$ does not start with $\tt AN$ (otherwise, $A=J_n$, which is singular).
By the Inverse Theorem,
$\epr(A^{-1})$ starts with $\tt ANSN$ and ends with $\tt A$.
It follows from Theorem \ref{ANSN} that either
$\epr(A^{-1}) = \tt ANSN\OL{SN}A$ or 
$\epr(A^{-1})=\tt ANSNAA$.
If it was the case that $\epr(A^{-1})=\tt ANSNAA$,
then the Inverse Theorem would imply that 
$A$ starts with $\tt AN$,
which would be a contradiction.
Thus, $\epr(A^{-1}) = \tt ANSN\OL{SN}A$,
implying that $n$ is odd, as desired.
By the Inverse Theorem, $\epr(A) = \tt \OL{NS}NSNAA$.
It follows from Theorem \ref{bipartite-cycle result} that 
$G(A) = C_n$.
Then, as $A$ is a 0--1 matrix, 
$A$ is permutationally similar to $A(C_n)$, as desired.

The last two statements of the theorem follow immediately from what was established above. 
\epf

Observe that Theorem \ref{NSNAA-SSNAA} implies that
the epr-sequences of symmetric 0--1 matrices of
order at least $5$ that end with $\tt NAA$ may be
partitioned into two classes:
those that end with $\tt NSNAA$ and
those that end with $\tt SSNAA$.
For any positive integer $k \geq 2$, $A(C_{2k+1})$ is
an example of a symmetric 0--1 matrix that
ends with $\tt NSNAA$ (see Theorem \ref{bipartite-cycle result}).
An example of a symmetric 0--1 matrix that
ends with $\tt SSNAA$, which we borrowed from \cite{BBH},
is the $9 \times 9$ matrix $A(C_3) \otimes A(C_{3})$,
whose epr-sequence is $\tt NSSSSSNAA$.

We are ready to establish the main result of this section,
which asserts that Conjecture \ref{NAA Conjecture} holds
if the matrix $A$ does not have both a zero and a nonzero principal minor of order $n-4$ (if $n \geq 5$).

\begin{thm}\label{main-result}
Let $n \geq 5$ be an integer,
let $A \in \zonn$ be symmetric and nonsingular, and
let $\epr(A) = \ell_1\ell_2 \cdots \ell_n$.
Suppose that $\ell_{n-4} \neq \tt S$.
Then the following statements are equivalent:
\ben
\item[$(i)$] $A$ ends with $\tt NAA$.
\vspace{-3mm}
\item[$(ii)$] $A^{-1}$ is equimodular with constant diagonal.
\een
\end{thm}

\bpf
By Remark \ref{Easy direction},
it suffices to show that (i) implies (ii).
Suppose that $A$ ends with $\tt NAA$.
By Remark \ref{(n-1)-minors}, it suffices to show that
all of the principal minors of order $n-1$ of $A$ are equal.
By Theorem \ref{NSNAA-SSNAA},
$A$ ends with $\tt NSNAA$ and
$A$ is permutationally similar to $A(C_n)$.
Then, as any two induced subgraphs of
$C_n$ of order $n-1$ are isomorphic,
and because $A$ is a 0--1 matrix, 
any two principal submatrices of $A$ of order $n-1$ are permutationally similar, implying that
any two principal minors of order $n-1$ of $A$ are equal.
\epf

With Theorem \ref{main-result} uncovered,
Conjecture \ref{NAA Conjecture} remains unresolved only in
the case where the matrix $A$ has the property of having
both a  zero and a nonzero principal minor of order $n-4$.

We shall now focus on
the determinant and order-$(n-1)$ principal minors of
symmetric matrices in $\zonn$ that end with $\tt NAA$,
particularly on their parity;
this is for two reasons.
In \cite[p. 55]{BBH}, the following question was formulated:
If $A \in \zonn$ is a nonsingular symmetric matrix such that
$A^{-1} \in \alphann$
(i.e., it is equimodular) with constant diagonal,
then what are the values that $\alpha$ may attain?
By Theorem \ref{adjoint-result},
$\alpha$ is either equal to, or the negative of, the ratio of
an order-$(n-1)$ principal minor of $A$ to $\det(A)$;
this gives us the first reason.
To see the second reason,
recall that Remark \ref{(n-1)-minors} implies that,
to establish Conjecture \ref{NAA Conjecture} (if it holds),
it suffices to establish the following statement:
If $A$ is a symmetric matrix in $\zonn$ that ends with $\tt NAA$,
then all of the order-$(n-1)$ principal minors of $A$ are equal.

To study the parity of the principal minors of
symmetric 0--1 matrices that end with $\tt NAA$,
we draw upon the author's work in \cite{XMR-Char 2}.
In Section \ref{s:intro}, we defined the
epr-sequence of a symmetric matrix over $\R$.
We now extend that definition to any field $\F$,
with the sole purpose of studying a (\textit{real}) 0--1 matrix by
interpreting it as a matrix over $\ZZ$.
If $A \in \zonn$ ($\subseteq \Rnn$) is symmetric,
then we shall denote by $\epr_2(A)$ its
epr-sequence as a matrix over $\ZZ$.
For example, if
\[
A=
\begin{pmatrix*}[r]
0 & 1 & 1 \\
1 & 0 & 1 \\
1 & 1 & 0
\end{pmatrix*} \in \{0,1\}^{3 \times 3} 
(\subseteq \R^{3 \times 3}),
\]
then $\epr(A) = \tt NAA$ and $\epr_2(A)=\tt NAN$.
Observe that, if $A$ is an arbitrary 
symmetric matrix in $\zonn$ ($\subseteq \Rnn$) and 
the $k$th term of $\epr_2(A)$ is $\tt A$ (respectively, $\tt N$),
then all of the principal minors of $A$ of order $k$ are odd (respectively, even);
and if the $k$th term is $\tt S$, then $A$ has both an odd and an even principal minor of order $k$.
\cite[Theorems 3.2, 3.8, 3.11]{XMR-Char 2} provide a
complete characterization of the epr-sequences of
symmetric matrices over $\ZZ$;
we draw upon these theorems to establish our next result.


\begin{prop}\label{epr2}
Let $A \in \zonn$ be symmetric and
$\epr(A) = \ell_1\ell_2 \cdots \ell_n$.
Suppose that $\ell_{n-2}\ell_{n-1}\ell_n = \tt NAA$.
Then $\epr_2(A)$ is of one of the following forms:
\[
\tt \OL{NS}NAN, \quad
\tt \OL{NS}NNN\OL{N}, \quad
\tt A\OL{S}NNN\OL{N}, \quad \text{or} \quad
\tt S\OL{S}NNN\OL{N}.
\]
\end{prop}

\bpf
Let $\epr_2(A) = \ell'_1\ell'_2\cdots \ell'_n$.
Observe that $\ell'_1 = \ell_1$.
Moreover, for all $j \in [n]$,
if $\ell_j=\tt N$, then $\ell'_j= \tt N$.
Thus, $\ell'_{n-2} = \tt N$.
We divide the proof into three cases, based on $\ell'_1$.

\noindent
\textbf{Case 1}: $\ell'_1 = \tt N$.

\noindent
As $\ell'_{n-2} = \tt N$,
it follows from \cite[Theorem 3.2]{XMR-Char 2}
that $\epr_2(A)$ is of one of the following forms:
$\tt NA\OL{NA}N$,
$\tt \OL{NS}NNN\OL{N}$ or
$\tt NS\OL{NS}NAN$.
If $\epr_2(A) = \tt NA\OL{NA}N$, then $A=J_n-I_n$,
implying that $\epr(A) = \tt NAA\OL{A}$
(which is not hard to verify) and, therefore, that $n=3$,
which implies that $\epr_2(A) = \tt NAN$.
Thus, $\epr_2(A)$ is of one of the following forms:
$\tt \OL{NS}NAN$ or $\tt \OL{NS}NNN\OL{N}$.

\noindent
\textbf{Case 2}: $\ell'_1 \neq \tt N$.

\noindent
As  $\ell'_{n-2} = \tt N$, it follows from 
\cite[Theorem 3.8 and Theorem 3.11]{XMR-Char 2}
that $\epr_2(A)$ is of one of the following forms:
$\tt A\OL{S}NNN\OL{N}$ or  $\tt S\OL{S}NNN\OL{N}$.
\epf

Each of the next two results is a corollary to
the result preceding it (they follow readily).

\begin{cor}
Let $A \in \zonn$ be symmetric and
$\epr(A) = \ell_1\ell_2 \cdots\ell_n$.
Suppose that $\epr(A)$ ends with $\tt NAA$.
Then $\epr_2(A)$ ends with either $\tt NAN$ or $\tt NNN$.
Moreover, if $n$ is even, then
$\epr_2(A)$ ends with $\tt NNN$.
Moreover, if $\ell_1 \neq \tt N$, then
$\epr_2(A)$ ends with $\tt NNN$.
\end{cor}

\begin{cor}\label{parity}
Let $A \in \zonn$ be symmetric and
$\epr(A) = \ell_1\ell_2 \cdots\ell_n$.
Suppose that $A$ ends with $\tt NAA$.
Then the following statements hold:
\ben
\item[$(i)$]
$\det(A)$ is even.
\item[$(ii)$]
If $n$ is even,
then all of the principal minors of $A$ of order $n-1$ are even.
\item[$(iii)$]
If $\ell_1 \neq \tt N$,
then all of the principal minors of $A$ of order $n-1$ are even.
\een
\end{cor}

As noted in \cite[p. 53]{BBH}, empirical evidence suggests that
nonsingular symmetric matrices in $\zonn$ whose
inverses are equimodular with constant diagonal
(which, by Remark \ref{Easy direction},
are matrices that end with $\tt NAA$) are rare.
The fact that the product of any
two order-$(n-1)$ principal minors
of any such matrix is a perfect square
(which is a fairly restrictive condition) is
supportive of such a hypothesis:

\begin{prop}\label{perfect square}
Let $A \in \zonn$ be symmetric.
Suppose that $A$ ends with $\tt NAA$.
Then the product of any two principal minors of $A$ of
order $n-1$ is a perfect square.
\end{prop}

\bpf
Let $i,j \in [n]$.
Let $d_i = \det(A(\{i\}))$, $d_j = \det(A(\{j\}))$ and
$d_{ij} = \det(A(\{i\},\{j\}))$.
We now show that $d_id_j = d_{ij}^{ \ 2}$.
As the desired conclusion follows immediately if $i=j$,
we assume that $i \neq j$.
Observe that $A$ is nonsingular (because it ends with $\tt A$).
Let $A^{-1} = [\alpha_{ij}]$.
By the Inverse Theorem, $A^{-1}$ starts with $\tt AN$.
Thus,
$0=A^{-1}[\{i,j\}] = \alpha_{ii}\alpha_{jj} -\alpha_{ij}^2$,
implying that $\alpha_{ii}\alpha_{jj} =\alpha_{ij}^2$.
As $A$ is symmetric, 
$d_{ij} = \det(A(\{j\},\{i\}))$.
It follows from Theorem \ref{adjoint-result}
that  $d_id_j = d_{ij}^{ \ 2}$.
\epf

Observe that Proposition \ref{perfect square} holds regardless of 
whether or not $A$ is a 0--1 matrix.

\section{0--1 matrices whose inverses are equimodular}\label{s:equim}
$\null$
\indent
This section focuses on 
nonsingular 0--1 matrices whose inverses are equimodular,
due to their relation to the second of the two families of 
matrices associated with Conjecture \ref{NAA Conjecture} 
(namely, the family of nonsingular symmetric 0--1 matrices
whose inverses are equimodular with constant diagonal).
In Sections \ref{s:AN and N} and \ref{s:NAA},
all of the matrices involved were symmetric,
as we are motivated by Conjecture \ref{NAA Conjecture}.
However, in this section, 
in Subsection \ref{ss:non-symmetric},
we shall start by establishing results that apply not only to symmetric matrices but also to non-symmetric matrices,
and then, in Subsection \ref{ss:symmetric and equim}, 
we shall revert our attention to symmetric matrices.
The present section will conclude with 
summative remarks about the section, 
in Subsection \ref{ss:summative}, 
by looking at its results through the lens of
Conjecture \ref{NAA Conjecture}.

\subsection{Not necessarily symmetric matrices}\label{ss:non-symmetric}
$\null$
\indent
The next fact, 
which generalizes \cite[Proposition 5]{BBH},
will help us cast some light upon
the structure of nonsingular 0--1 matrices
whose inverses are equimodular,
and its proof, which we include merely for completeness,
is the same as that of \cite[Proposition 5]{BBH}.

\begin{prop}\label{even-row-sums}
Let $n \geq 2$ be an integer and $A \in \Rnn$ be nonsingular.
Suppose that each nonzero entry of $A$ is an odd integer.
Suppose that $A^{-1} = qN$, for
some rational number $q$ and
some matrix $N$ all of whose entries are odd.
Then each row, and each column, of $A$ has
an even number of nonzero entries.
\end{prop}

\bpf
Let $A= [\vec{a}_1, \vec{a}_2, \dots, \vec{a}_n]$ and 
$N^T= [\vec{w}_1, \vec{w}_2, \dots, \vec{w}_n]$.
Let $j \in [n]$ and $m$ be 
the number of nonzero entries of $\vec{a}_j$.
Let $i \in [n] \setminus \{j\}$.
Observe that 
$\vec{w}_i^T\vec{a}_j$ is the $(i,j)$-entry of $NA$.
Then, as $NA$ is a diagonal matrix,
and because $i \neq j$, $\vec{w}_i^T\vec{a}_j = 0$, 
implying that $\vec{w}_i^T\vec{a}_j$ is even.
As $\vec{w}_i^T\vec{a}_j$ is
the sum of $m$ odd numbers,
the parity of $\vec{w}_i^T\vec{a}_j$ is
the same as that of $m$.
It follows that $m$ is even.
Thus, each column of $A$ has 
an even number of nonzero entries.

Applying what was established in the previous paragraph to
$A^T$ shows that each column of 
$A^T$ has an even number of nonzero entries, implying that
each row of $A$ has an even number of nonzero entries.
\epf

Each row and column of a 
0--1 matrix (of order at least $2$) whose inverse is equimodular 
has an even number of nonzero entries:

\begin{cor}\label{Rows-and-columns-if-inverse-is-equim}
Let $n \geq 2$ be an integer and $A \in \zonn$ be nonsingular.
Suppose that $A^{-1}$ is equimodular.
Then each row, and each column, of $A$ has
an even number of nonzero entries.
\end{cor}

\bpf
It is clear that $A^{-1}$ being equimodular implies that
$A^{-1}=\alpha N$, for some number $\alpha$ and
some matrix $N \in \{-1,1\}^{n \times n}$.
Observe that each entry of $N$ is odd.
By Theorem \ref{adjoint-result},
the fact that each entry of $A$ is rational implies that
each entry of $A^{-1}$ is as well,
implying that $\alpha$ is rational.
It follows from Proposition \ref{even-row-sums} that 
each row, and each column, of $A$ has 
an even number of nonzero entries.
\epf

The determinants of nonsingular 0--1 matrices whose 
inverses are equimodular are all even:

\begin{thm}\label{1/even}
Let $n \geq 2$ be an integer and $A \in \zonn$ be nonsingular.
Suppose that $A^{-1}$ is equimodular,
with $A^{-1} \in \alphann$, for some $\alpha$.
Then the following statements hold:
\ben
\item\label{1/2m}
$\alpha = \frac{1}{2m}$, for some integer $m$.
\item\label{order n-1 divides det}
Each minor of $A$ of order $n-1$ divides $\det(A)$.
\item\label{det is even}
$\det(A)$ is even.
\een
\end{thm}

\bpf
As $A^{-1} \in \alphann$, 
there exists a matrix $N \in \{-1,1\}^{n \times n}$ such that
$A^{-1} = \alpha N$.
As $A^{-1}$ is nonsingular, $\alpha \neq 0$.
Thus, $AN = \frac{1}{\alpha} I_n$.
By Corollary \ref{Rows-and-columns-if-inverse-is-equim},
each row of $A$ has an even number of nonzero entries.
Then, as $A \in \zonn$ and $N \in \{-1,1\}^{n \times n}$,
each entry of $AN$ is 
the sum of an even number of odd numbers;
as such a sum is even, each entry of $AN$ is even.
Then, as $AN = \frac{1}{\alpha} I_n$,
$\frac{1}{\alpha}$ is an even integer.
Hence, $\alpha = \frac{1}{2m}$, for some integer $m$,
establishing statement (\ref{1/2m}).

Let $d = \det(A)$ and $t$ be a minor of $A$ of order $n-1$.
As $A^{-1} \in \alphann$,
Theorem \ref{adjoint-result} implies that
$\big|\frac{t}{d}\big| = |\alpha|$.
Thus, $\big|\frac{d}{t}\big| = |2m|$.
Then, as $d$ and $t$ are integers (because $A$ is a 0--1 matrix),
$t$ divides $d$. 
Hence, statement (\ref{order n-1 divides det}) follows.
As $|d| = |2mt|$, and because $d$ and $t$ are integers, statement (\ref{det is even}) holds.
\epf

We note that Theorem \ref{1/even} asserts that
not only 
all principal minors of $A$ of order $n-1$ divide $\det(A)$, 
but also 
all non-principal minors of $A$ of order $n-1$ divide $\det(A)$ .

Due to the seeming irrelevance of the next fact, 
it is worthwhile to note that its inclusion will be justified at 
the end of this section (in Remark \ref{small order remark}).

\begin{prop}\label{det=2k}
Let $n \geq 2$ be an integer and $A \in \zonn$ be nonsingular.
Suppose that $A^{-1}$ is equimodular,
with $A^{-1} \in \alphann$, for some $\alpha$.
Suppose that $|\det(A)| = 2^k$, for some integer $k >0$.
Then there exists an integer $t$ with $0 \leq t < k$ such that 
the modulus of each minor of $A$ of order $n-1$ is $2^t$.
Moreover, $|\alpha| = \frac{1}{2^{k-t}}$.
Moreover, if $A$ has an odd minor of order $n-1$,
then $t=0$ and $|\alpha| = \big|\frac{1}{\det(A)}\big|$.
\end{prop}

\bpf
As $A^{-1}$ is equimodular,
Theorem \ref{adjoint-result} implies that the moduli of 
any two minors of $A$ of order $n-1$ are equal.
By Theorem \ref{1/even},
each minor of $A$ of order $n-1$ divides $\det(A)$.
Thus, there exists an integer $t$ with $0 \leq t \leq k$ such that 
the modulus of each minor of $A$ of order $n-1$ is $2^t$.
By Theorem \ref{adjoint-result}, 
$|\alpha| = \frac{2^t}{2^{k}} = \frac{1}{2^{k-t}}$.
By Theorem \ref{1/even},
$|\alpha|$ is the reciprocal of an even integer.
Thus, $t \neq k$, implying that $0 \leq t <k$, as desired. 
The last statement follows readily from 
the statements preceding it.
\epf

The next result asserts that the inverse of a 
Schur complement in a nonsingular  0--1 matrix whose
inverse is equimodular is itself equimodular,
and it follows immediately from 
the Schur Complement Theorem.

\begin{prop}\label{Schur-with-equim}
Let $n\geq 2$ be an integer and $A \in \zonn$ be nonsingular, and let $\mu \subset [n]$.
Suppose that $A[\mu]$ is nonsingular, and
let $C=A/A[\mu]$.
If $A^{-1}$ is equimodular, then $C^{-1}$ is equimodular.
Furthermore, if
$A^{-1} \in \alphann$, then
$C^{-1} \in \{-\alpha,\alpha\}^{(n-|\mu|)\times(n-|\mu|)}$.
\end{prop}

We note that the statement of 
Proposition \ref{Schur-with-equim} would 
continue to hold if $A$ was replaced with an
arbitrary nonsingular matrix in $\Rnn$.

\begin{lem}\label{Z-xy for column difference}
Let $n \geq 3$ be an integer, and
let $A \in \zonn$ be nonsingular and of the form
\begin{equation*}
A=
\left(
\begin{array}{c|c}
 1 & \vec{y}^T  \\
 \hline
 \vec{x} & Z
\end{array}
\right),
\end{equation*}
for some
$\vec{x}, \vec{y}  \in \R^{n-1}$
and \
$Z \in \R^{(n-1) \times (n-1)}$.
Suppose that $A^{-1}$ is equimodular.
Then each row, and each column, of $Z-\vec{x}\vec{y}^T$ has
an even number of nonzero entries.
\end{lem}

\bpf
Let $C=Z-\vec{x}\vec{y}^T$ and $\mu = \{1\}$.
Observe that $Z = A[\mu^c]$ 
(where $\mu^c = [n] \setminus \mu$), 
that $\vec{x} = A[\mu^c, \mu]$,
that $\vec{y}^T = A[\mu, \mu^c]$,
that $A[\mu]$ is the $1 \times 1$ matrix $\mtx{1}$,
which is nonsingular, and that $(A[\mu])^{-1} =\mtx{1}$.
It follows, then, that $C = A/A[\mu]$.
As $A^{-1}$ is equimodular, 
Proposition \ref{Schur-with-equim} implies that
$C^{-1}$ is equimodular as well,
implying that $C^{-1} = \alpha N$, for some number
$\alpha$ and some $N \in \{-1,1\}^{(n-1)\times(n-1)}$.
As each entry of $C$ is rational,
Theorem \ref{adjoint-result} implies that
each entry of $C^{-1}$ is also rational and, therefore,
that $\alpha$ is rational.
As $Z$ and $\vec{x}\vec{y}^T$ are both 0--1 matrices, 
$C \in \{-1,0,1\}^{(n-1) \times (n-1)}$.
Then, as $n-1 \geq 2$,
Proposition \ref{even-row-sums} implies that
each row, and each column, of $C$ has 
an even number of nonzero entries.
\epf

Our next result casts light upon the structure of
nonsingular 0--1 matrices whose inverses are equimodular.

\begin{thm}\label{column difference}
Let $n \geq 3$ be an integer and $A \in \zonn$ be nonsingular, and let
$A=[\vec{a}_1, \vec{a}_2, \dots, \vec{a}_n]$ and
$A^T=[\vec{r}_1, \vec{r}_2, \dots, \vec{r}_n]$.
Suppose that $A^{-1}$ is equimodular.
Then each of
$\vec{a}_i - \vec{a}_j$  and
$\vec{r}_i - \vec{r}_j$
has an even number of nonzero entries, for all $i,j \in [n]$.
\end{thm}

\bpf
Let $A=[a_{ij}]$.
We start by showing that 
$\vec{a}_i - \vec{a}_j$ has an even number of nonzero entries,
for all $i,j \in [n]$.
Let $i,j \in [n]$.
As the desired conclusion follows immediately if $i=j$,
assume that $i \neq j$.
We consider two cases, based on $\vec{a}_i^T\vec{a}_j$.

\noindent
\textbf{Case 1}: $\vec{a}_i^T\vec{a}_j = 0$.

\noindent
Let $m$ and $\ell$ be the number of nonzero entries of 
$\vec{a}_i$ and $\vec{a}_j$, respectively.
As $\vec{a}_i^T\vec{a}_j = 0$,
and because $\vec{a}_i$ and $\vec{a}_j$ are 0--1 vectors, 
the number of nonzero entries of 
$\vec{a}_i - \vec{a}_j$ is $m + \ell$.
By Corollary \ref{Rows-and-columns-if-inverse-is-equim}, 
$m$ and $\ell$ are both even, implying that 
$m + \ell$ is even, as desired.

\noindent
\textbf{Case 2}: $\vec{a}_i^T\vec{a}_j \neq 0$.

\noindent
As $\vec{a}_i$ and $\vec{a}_j$ are 0--1 vectors,
there exists $k \in [n]$ such that $a_{ki} = a_{kj} = 1$
(otherwise, $\vec{a}_i^T\vec{a}_j = 0$).
Let $P \in \Rnn$ be a permutation matrix such that 
the first row of $PA$ is the $k$th row of $A$.
Let $Q \in \Rnn$ be a permutation matrix such that 
the first and second columns of $PAQ$ are, respectively,
the $i$th and $j$th columns of $PA$.
Let $B = PAQ$, 
$B=[\vec{b}_1,\vec{b}_2,\dots, \vec{b}_n]$ and $B=[b_{ij}]$. 
Observe that the number of nonzero entries of 
$\vec{a}_i - \vec{a}_j$ 
is equal to the number of nonzero entries of 
$\vec{b}_1 - \vec{b}_2$.
Thus, it suffices to show that $\vec{b}_1 - \vec{b}_2$ has 
an even number of nonzero entries.
Observe that 
$b_{11} = 1$ (because $a_{ki} = 1$) and
$b_{12} = 1$ (because $a_{kj} = 1$).
Thus, $B$ is of the form
\begin{equation*}
B=
\left(
\begin{array}{c|c}
 1 & \vec{y}^T  \\
 \hline
 \vec{x} & Z
\end{array}
\right),
\end{equation*}
for some $\vec{x}, \vec{y} \in \R^{n-1}$ 
and \ $Z \in \R^{(n-1) \times (n-1)}$.
Let
$\vec{y} = [y_1, \dots, y_{n-1}]^T$  and 
$Z = [\vec{z}_1, \dots, \vec{z}_{n-1}]$.
Then, as $y_1 = b_{12} = 1$,
\begin{equation*}
\vec{b}_1 - \vec{b}_2=
\left(
\begin{array}{c}
 1 - y_1  \\
 \hline
 \vec{x}  - \vec{z}_1 
\end{array}
\right) 
=
-\left(
\begin{array}{c}
 0  \\
 \hline
 \vec{z}_1 - \vec{x} 
\end{array}
\right).
\end{equation*}
Thus, it suffices to show that $\vec{z}_1 - \vec{x} $ has 
an even number of nonzero entries.
Observe that $B^{-1} = Q^TA^{-1}P^T$.
It follows, then, that $B^{-1}$ is equimodular 
(because $A^{-1}$ is equimodular).
Then, as $B \in \zonn$, and because $n \geq 3$,
Lemma \ref{Z-xy for column difference} implies that
$\vec{z}_1 - y_1\vec{x}$, 
the first column of $Z - \vec{x}\vec{y}^T$,
has an even number of nonzero entries;
as $y_1 = 1$, $\vec{z}_1 - \vec{x}$ has an 
even number of nonzero entries, as desired.

Let $M=A^T$. 
Observe that $M \in \zonn$, that $M$ is nonsingular, 
that $M = [\vec{r}_1, \vec{r}_2, \dots, \vec{r}_n]$,
and that 
$M^{-1}$ is equimodular ($M^{-1}=(A^{-1})^T$).
It follows from what was established in 
the previous paragraph that, for all $i,j \in [n]$,
$\vec{r}_i - \vec{r}_j$ has an even number of nonzero entries.
\epf

\subsection{Symmetric matrices}\label{ss:symmetric and equim}
$\null$
\indent
We shall now revert our attention to symmetric matrices,
albeit with our focus still on 0--1 matrices whose inverses are equimodular.

\begin{lem}\label{even number diag lemma}
Let $B \in \Rnn$ be symmetric.
Suppose that each row of $B$ has 
an even number of nonzero entries.
Then $B$ has an even number of nonzero diagonal entries.
\end{lem}

\bpf
Let $m$ be the number of nonzero diagonal entries of $B$.
Let $M$ be the matrix that results from
replacing each nonzero diagonal entry of $B$ with $0$.
Then, as each row of $B$ has an 
even number of nonzero entries,
$M$ has exactly $m$ rows with
an odd number of nonzero entries.
It follows, then,
that $G(M)$ has exactly $m$ vertices of odd degree.
By the (well-known) degree-sum formula, $m$ is even.
\epf


A nonsingular symmetric 0--1 matrix (of order at least two) whose inverse is equimodular has
an even number of nonzero diagonal entries:

\begin{thm}\label{Diagonal-if-inverse-is-equim}
Let $n \geq 2$ be an integer, and
let $A \in \zonn$ be symmetric and nonsingular.
Suppose that $A^{-1}$ is equimodular.
Then $A$ has an even number of nonzero diagonal entries.
\end{thm}

\bpf
By Corollary \ref{Rows-and-columns-if-inverse-is-equim},
each row of $A$ has an even number of nonzero entries.
The desired conclusion follows from 
Lemma \ref{even number diag lemma}.
\epf

Lifting the assumption that $A$ is symmetric in 
Theorem \ref{Diagonal-if-inverse-is-equim} renders 
the statement of the proposition false:

\begin{ex}\label{diag without symmetry example}
\rm
Observe that the inverse of the non-symmetric 0--1 matrix
\begin{equation*}
A=
\left(
\begin{array}{ccc}
 1 & 1 & 0 \\
 0 & 1 & 1 \\
 1 & 0 & 1 \\
\end{array}
\right) 
\text{\quad is \quad}
A^{-1} = 
\frac{1}{2}
\left(
\begin{array}{ccc}
 1 & -1 & 1 \\
 1 & 1 & -1 \\
 -1 & 1 & 1 \\
\end{array}
\right),
\end{equation*}
which is equimodular.
As $A$ has an odd number of nonzero diagonal entries,
we conclude that a non-symmetric 0--1 matrix whose 
inverse is equimodular need not have an 
even number of nonzero diagonal entries.
\end{ex}

The inclusion of the following observation will be justified shortly.

\begin{obs}\label{obs for A-aa}
Let $A \in \zonn$ be symmetric, $A=[a_{ij}]$ and 
$A=[\vec{a}_1, \vec{a}_2, \dots, \vec{a}_n]$.
Let $k \in [n]$, $B = A-\vec{a}_k\vec{a}_k^T$ and
$B=[\vec{b}_1, \vec{b}_2, \dots, \vec{b}_n]$.
Then, for all $i \in [n]$, 
\begin{equation*}
\vec{b}_i = 
\begin{cases}
\vec{a}_i - \vec{a}_k,  & \text{if \ $a_{ik} = 1$;} \\ 
\vec{a}_i,                   & \text{if \ $a_{ik} = 0$}.
\end{cases}
\end{equation*}
\end{obs}

The purpose of making the observation above 
(Observation \ref{obs for A-aa}) is to note that
knowledge of the structure of $A$ may 
be revealed through the structure of $A-\vec{a}_k\vec{a}_k^T$,
as the matrix $A-\vec{a}_k\vec{a}_k^T$ is obtained from 
$A$ by simply subtracting $\vec{a}_k$ from 
certain column vectors of $A$;
this observation is rendered useful, for us, by 
the following theorem,
which should be viewed through the lens of 
Observation \ref{obs for A-aa}.

\begin{thm}\label{A-aa}
Let $n \geq 3$ be an integer, 
let $A \in \zonn$ be symmetric and nonsingular, and 
let $A=[\vec{a}_1, \vec{a}_2, \dots, \vec{a}_n]$.
Suppose that $A^{-1}$ is equimodular. 
Then $A-\vec{a}_k\vec{a}_k^T$ has 
an even number of nonzero diagonal entries,
for all $k \in [n]$.
\end{thm}

\bpf
Let $k \in [n]$, 
$B=A-\vec{a}_k\vec{a}_k^T$ and 
$B=[\vec{b}_1, \vec{b}_2, \dots, \vec{b}_n]$.
It follows from Observation \ref{obs for A-aa} that, 
for all $i \in [n]$,
either
$\vec{b}_i = \vec{a}_i - \vec{a}_k$ or
$\vec{b}_i = \vec{a}_i$.
Thus, by
Corollary \ref{Rows-and-columns-if-inverse-is-equim} and
Theorem \ref{column difference},
each column of $B$ has an even number of nonzero entries.
Then, as $B$ is symmetric, 
each row of $B$ has an even number of nonzero entries.
By Lemma \ref{even number diag lemma},
$B$ has an even number of nonzero diagonal entries.
\epf

Unsurprisingly, the assumption that $A$ is symmetric in 
Theorem \ref{A-aa} is important:

\begin{ex}
\rm
Let $A$ be the $3 \times 3$ matrix in
Example \ref{diag without symmetry example},
and let $A=[\vec{a}_1,\vec{a}_2, \vec{a}_3]$.
Observe that $A$ is a non-symmetric 0--1 matrix whose 
inverse is equimodular, and that 
\begin{equation*}
A - \vec{a}_1\vec{a}_1^T=
\left(
\begin{array}{ccc}
 0 & 1 & -1 \\
 0 & 1 & 1 \\
 0 & 0 & 0 \\
\end{array}
\right).
\end{equation*}
Thus, $A - \vec{a}_1\vec{a}_1^T$ has an 
odd number of nonzero diagonal entries.
We conclude that relaxing the hypothesis of 
Theorem \ref{A-aa} by not requiring the matrix $A$ to 
be symmetric leads to a false statement.
\end{ex}

\subsection{Summative remarks}\label{ss:summative}
$\null$
\indent
We shall now provide some summative remarks
about this section (Section \ref{s:equim}).
It is clear that the main results of this section apply to 
the matrices in the family of 
nonsingular symmetric 0--1 matrices whose 
inverses are equimodular with constant diagonal.
Although the results of this section apply to 
other matrices besides those in this family, 
we shall keep this family at the center of 
our summative remarks below,
due to its relation to Conjecture \ref{NAA Conjecture},
which was our primary motivation for investigating
(not necessarily symmetric) 0--1 matrices whose 
inverses are equimodular 
(with not necessarily constant diagonal) in this section.

Hitherto, 
what appears on the literature regarding matrices in 
the family of nonsingular symmetric 0--1 matrices whose 
inverses are equimodular with constant diagonal is 
largely focused on those that have zero diagonal
(which may be interpreted as
adjacency matrices of graphs without loops).
On the other hand, of those that have 
at least one nonzero diagonal entry (examples of which are
the last three matrices in \cite[p. 54]{BBH}),
not much is known, apart from what is implied by
\cite[Proposition 5]{BBH}
(i.e., that each of their rows and columns has 
an even number of nonzero entries).
It is, therefore, worthwhile to note that 
the results uncovered in this section apply to 
\textit{arbitrary} nonsingular symmetric 0--1 matrices whose 
inverses are equimodular with constant diagonal
(meaning that there is no restriction imposed on their diagonals);
of particular relevance are
Theorems \ref{column difference}, 
\ref{Diagonal-if-inverse-is-equim} and \ref{A-aa},
as they cast light upon the structure of these matrices
(Theorem \ref{A-aa} should be viewed through 
the lens of Observation \ref{obs for A-aa}).

The following remark,
whose verification is a trivial computational exercise
(albeit with the help of \cite[p. 53--54]{BBH}),
will help us elucidate the relevance of 
Proposition \ref{det=2k}, 
by linking it to the family of 
nonsingular symmetric 0--1 matrices whose 
inverses are equimodular with constant diagonal.

\begin{rem}\label{small order remark}
\rm
Let $3 \leq n \leq 9$ be an integer, and 
let $A \in \zonn$ be symmetric and nonsingular.
Suppose that $A^{-1}$ is equimodular with constant diagonal.
Then $|\det(A)| = 2^k$, for some integer $k>0$.
Moreover, if $A$ is neither 
permutationally similar to $A(C_3) \otimes A(C_3)$ nor 
permutationally similar to
one of the last two $9 \times 9$ matrices on \cite[p. 54]{BBH},
then $A$ has an odd minor of order $n-1$. 
\end{rem}

Up to permutation similarity,
there is a total of fifteen nonsingular symmetric 0--1 matrices 
of order at least $3$ and at most $9$ 
whose inverses are equimodular with constant diagonal 
\cite[p. 53--54]{BBH}
(each of these fifteen matrices appears in \cite{BBH}).
Remark \ref{small order remark} asserts that,
if $A \in \zonn$ is one of the aforementioned fifteen matrices
but is neither 
permutationally similar to $A(C_3) \otimes A(C_3)$ nor 
permutationally similar to
one of the last two $9 \times 9$ matrices on \cite[p. 54]{BBH},
then $A$ satisfies \textit{all} of 
the hypotheses of Proposition \ref{det=2k} and, therefore,
$|\alpha| = \big|\frac{1}{\det(A)}\big|$, 
where $A^{-1} \in \alphann$.

\end{document}